\documentclass{article}
\usepackage[utf8]{inputenc}
\usepackage{amsmath}
\usepackage{amsfonts}
\usepackage{mathrsfs}
\newtheorem{definition}{Definition}
\newtheorem{proposition}{Proposition}
\newtheorem{theorem}{Theorem}
\newtheorem{conjecture}{Conjecture}
\newtheorem{lemma}{Lemma}
\usepackage[margin=1in]{geometry}
\usepackage{graphicx}
\usepackage{url}
\usepackage{authblk}

\title{An uncountable union of line segments with null two-dimensional measure}
\author{Parker Kuklinski}
\date{October 2023}

\begin{document}

\maketitle

\begin{abstract}
In this paper we construct an uncountable union of line segments $T$ which has full intersection with the sets $\left(\{0\}\times [0,1]\right) \cup \left(\{1\}\times [0,1]\right) \subset\mathbb{R}^2$ but has null two-dimensional measure. Further results are proved on the decay rate of $\mu (T)$ if the line segments comprising $T$ are replaced with increasingly fine approximations by parallelograms.
\end{abstract}

\section{Introduction}

In \cite{kaczynski67}, Kaczynski proves a variety of results placing restrictions on possible boundary functions of a half-plane function with some property. The last result is different; he constructs a measurable half plane function $f$ with a nonmeasurable boundary function. Crucial to this construction is a measure theoretic generalization of the trapezoid formula which we document here. For the remainder of the paper, let $X_0=\{ (x,y):y=0\}$ and $X_1=\{ (x,y):y=1\}$, let $\mu$ refer to Lebesgue measure in the contextually obvious dimension (we will distinguish $\mu _1$ and $\mu _2$ to refer to one and two-dimensional Lebesgue measure respectively when necessary), and let $\mu^*$ refer to outer measure.
\begin{definition}
A \emph{trapezoid} $T=\cup\ell$ is a disjoint union of line segments with each segment having its two endpoints on $X_0$ and $X_1$ (i.e. $\ell =\{((1-t)x_0+tx_1,t):t\in [0,1]\}$).
\end{definition}
If $T$ is a trapezoid, let $T_0=T\cap X_0$ and $T_1=T\cap X_1$. In \cite{kaczynski67}, Kaczynski proved the following result:
\begin{theorem}
If $T$ is a trapezoid, then $\mu ^*_2(T)=(\mu ^*_1(T_0)+\mu ^*_1(T_1))/2$.
\end{theorem}
This theorem was used to construct a collection of non-intersecting paths (not lines) $\gamma _x$ corresponding to each point $x\in [0,1]$ such that for each path $\lim _{t\rightarrow 0} \gamma _x(t)=(x,0)$ but the union satisfies $\mu\left(\cup _{x\in [0,1]}\gamma_x\right) =0$. This construction immediately leads to the creation of a measurable half plane function with non-measurable boundary function \cite{kaczynski67}. We depict the situation in Figure 1.

\begin{figure}
\centering
\includegraphics[scale=0.25]{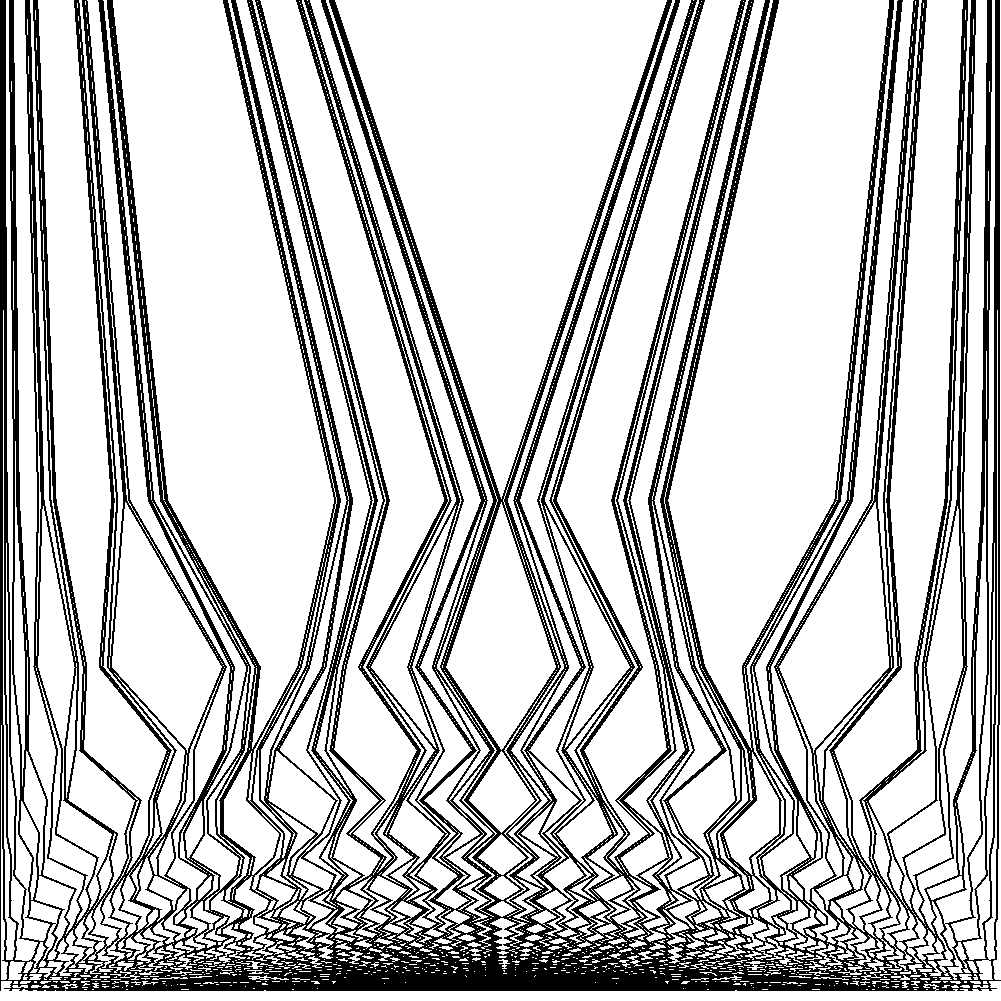}
\caption{Plot of $\cup \gamma _x$ for $x\in [0,1]$ where $\gamma _{x_1}\cap\gamma _{x_2}=\emptyset$, $\lim _{t\rightarrow 0} \gamma _x(t)=(x,0)$, but $\mu\left(\cup _{x\in [0,1]}\gamma_x\right) =0$. Given a nonmeasurable function $\varphi (x)$, if we define $f:\left(\mathbb{R}\times (0,\infty)\right)\rightarrow\mathbb{R}$ as $f(\gamma _x)=\varphi (x)$ and vanishing elsewhere, then $f$ is measurable in the upper half plane but has a nonmeasurable boundary function $\varphi$.}
\end{figure}

One may wonder the extent to which this result holds if we omit the restriction that the lines of a trapezoid are not disjoint. To this end, we amend the trapezoid definition:
\begin{definition}
An \emph{unrestricted trapezoid} $T=\cup\ell$ is a union (not necessarily disjoint) of line segments with each segment having its two endpoints on $X_0$ and $X_1$ (i.e. $\ell =\{((1-t)x_0+tx_1,t):t\in [0,1]\}$).
\end{definition}
In \cite{kaczynski67}, Kaczynski proves that there is no sensible upper bound in terms of $\mu _1(T_0)$ and $\mu _1 (T_1)$ on the two dimensional measure of an unrestricted trapezoid.
\begin{theorem}
There exists an unrestricted trapezoid $T$ with $\mu _1(T_0)=\mu _1(T_1)=0$ and $\mu _2(T)=\infty$.
\end{theorem}
{\bf Proof:} Let $M$ be a residual set of measure zero (e.x. let $x_k$ be an enumeration of rationals and $M=\bigcap _{n=1}^\infty\bigcup_{k=1}^\infty\left( x_k-\frac{1}{2^{k+1}n},x_k+\frac{1}{2^{k+1}n}\right)$). For any point $(x,y)$ ($y\in (0,1)$) and line $\ell$ passing through $(x,y)$ with angle $\theta$ with the $x$-axis, let $F_0(\theta )=\left( x-y\cot(\theta ),0\right)$ and $F_1(\theta )=\left( x+(1-y)\cot (\theta ),1\right)$ denote the intersections $\ell\cap X_0$ and $\ell\cap X_1$ respectively. $F_0$ and $F_1$ are homeomorphisms from $(0,\pi )$ onto $X_0$ and $X_1$ respectively, so $F_0^{-1}(M)$ and $F_1^{-1}(M)$ are both residual sets of measure zero. Since the intersection of residual sets is nonempty, for every $(x,y)$ with $y\in (0,1)$ we can choose an angle $\alpha =F_0^{-1}(M)\cap F_1^{-1}(M)$ such that the line $\ell$ intersecting $(x,y)$ with angle $\alpha$ intersects $X_0\cap M$ and $X_1\cap M$. Let the unrestricted trapezoid $T$ be the union of all such lines so that $\mu _1(T_0)=\mu _1(T_1)=\mu _1(M)=0$ and $\mu _2(T)=\mu _2(\{(x,y),y\in (0,1)\} )=\infty$. $\hfill\Box$.

\section{No lower bound}

In \cite{kaczynski67}, Kaczynski poses the opposite question as an open problem, asking if there is a lower bound on the two-dimensional measure of an unrestricted trapezoid in terms of $\mu (T_0)$ and $\mu (T_1)$. We will prove the following theorem indicating no lower bound exists:
\begin{theorem}
There exists an unrestricted trapezoid $T$ with $\mu _1(T_0)=\mu _1(T_1)=\infty$ and $\mu _2(T)=0$.
\end{theorem}
Since it suffices to find an unrestricted trapezoid with measure zero and sets $\mu (T_0)=\mu (T_1)=\mu ([0,1])=1$, we can formulate the question a different way:
\begin{proposition}
There exists a bijection $f:[0,1]\rightarrow [0,1]$ such that, for the set
$$S=\bigcup _{t\in [0,1]}\bigcup _{x\in [0,1]}\left( (1-t)x+tf(x),t\right),$$
we have $\mu (S)=0$.
\end{proposition}
It is clear that Theorem 3 is a corollary of Proposition 1.

It is perhaps not immediate that such a function should exist. Any monotone $f$ induces a standard trapezoid which abides by Theorem 1. The choice $f(x)=1-x$ gives measure $\mu (S)=1/2$, and other decreasing functions cannot be massaged to induce a smaller measure. We will instead draw inspiration from fractal geometry \cite{falconer97}; the particular iterative scheme we choose will create one-dimensional slices in $x$ similar to the Cantor middle third set. In Figure 2, the basic pattern is represented by the piecewise function
\[   f_1(x)=\left\{
\begin{array}{ll}
      x & x\in [0,\frac{1}{3}) \\
      x+\frac{1}{3} & x\in [\frac{1}{3},\frac{2}{3}) \\
      x-\frac{1}{3} & x\in [\frac{2}{3},1] 
\end{array} 
\right. \]
otherwise we divide $[0,1]$ into three subintervals and connect them via parallelograms where $f_1$ essentially swaps the second and the third. In this case, $\mu (S)=5/6$. In each subsequent iteration we embed an affine transformed version of the base pattern in each parallelogram, thus reducing the total measure of the unrestricted trapezoid (we will prove this limits to zero as in Proposition 1).

\begin{figure}
\centering
\includegraphics[scale=0.2]{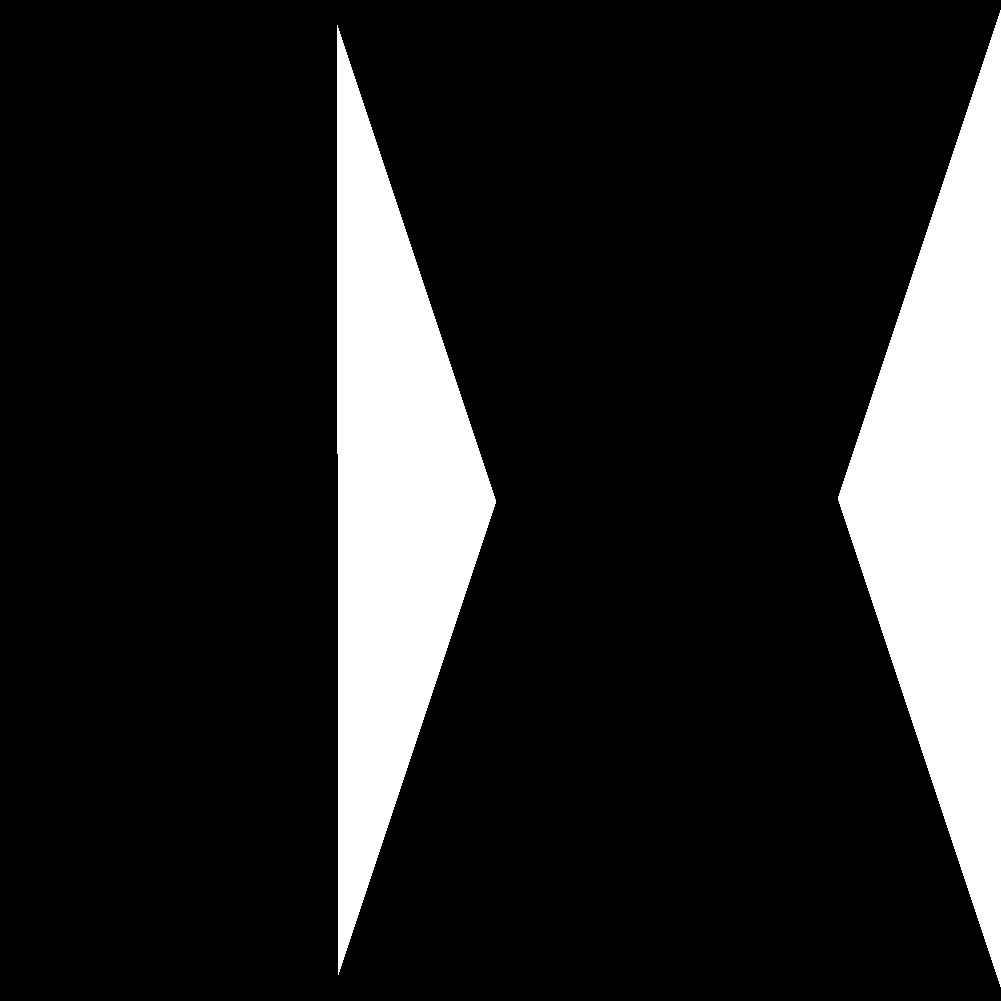}\hspace{1cm}\includegraphics[scale=0.2]{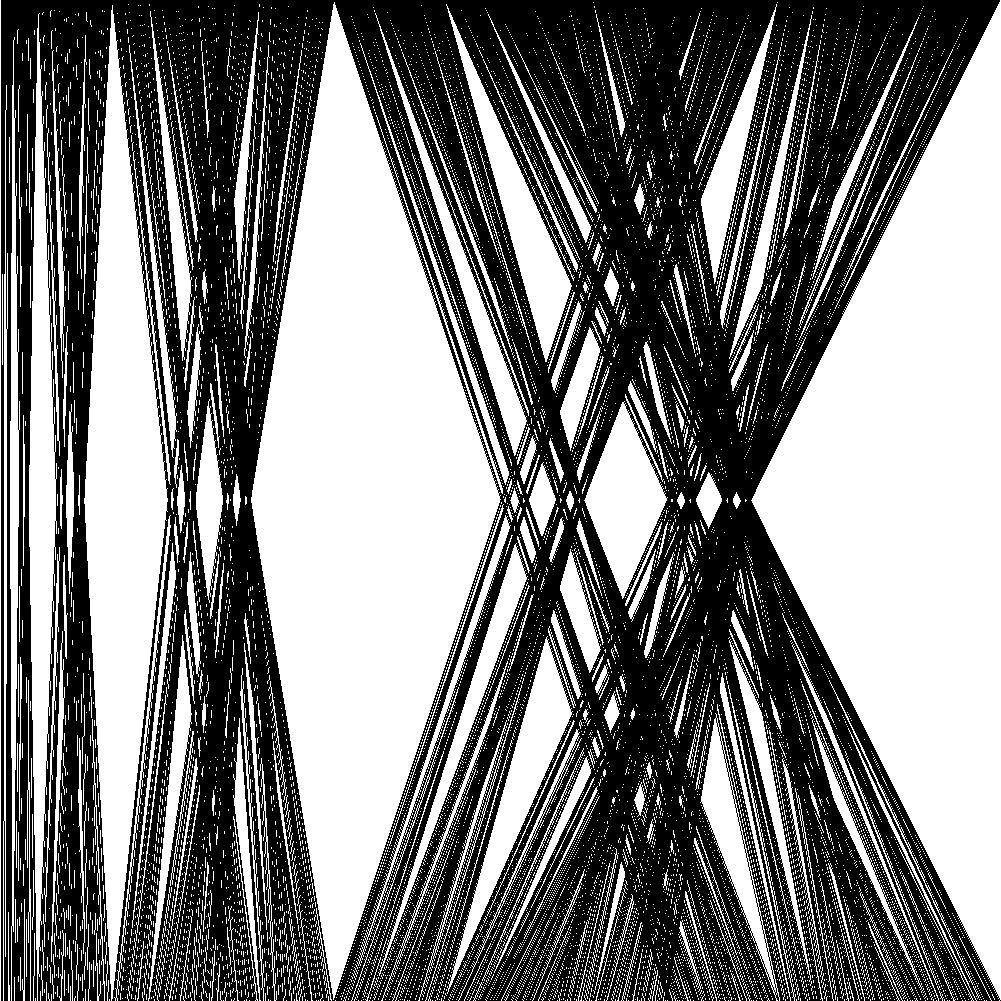}
\caption{Examples of the fractal construction of an unrestricted trapezoid with null measure. (\emph{Left}) The base pattern induced by $f_1$ (\emph{Right}) Illustration of unrestricted trapezoid induced by $f$.}
\end{figure}

Since this procedure iteratively divides subintervals of $[0,1]$ in 3, it is perhaps not surprising that this this sequence of functions $\{ f_k\}$ has a natural representation on base 3 fractions, namely
$$f_k(0.x_1x_2x_3,,,)=0.y_1y_2...y_kx_{k+1}x_{k+2}...$$
where $y_j=2x_j~(\text{mod }3)$. Let $f=\lim _{k\rightarrow\infty}f_k$ such that for input $x\in [0,1]$, $f$ replaces every 1 in its base 3 expansion with a 2 and vice-versa. We will prove that this function satisfies proposition 1.

Before we prove proposition 1, we will state a result from \cite{kenyon97} on modified Cantor sets:
\begin{proposition}
Let $C_t$ be a modified Cantor set with representation
$$C_t=\left\{\sum _{k=1}^\infty\frac{x_k}{3^k}:x\in\{ 0,1,t\}\right\} .$$
Then the following holds:
\[   \mu _1(C_t)=\left\{
\begin{array}{ll}
      1/q & t=p/q\text{ in lowest terms, }p+q\equiv 0\text{ (mod 3)} \\
      0 & \text{otherwise}
\end{array} 
\right. \]
\end{proposition}

We are now ready to prove Proposition 1.

\noindent {\bf Proof of Proposition 1:} Let us define $f$ acting on base 3 fractions as
$$f_k(0.x_1x_2x_3,,,)=0.y_1y_2y_3...$$
where $y_j=2x_j~(\text{mod }3)$. Then the linear combination $(1-t)x+tf(x)$ satisfies
$$S_t=\bigcup _{x\in [0,1]}((1-t)x+tf(x))=\left\{\sum _{k=1}^\infty\frac{x_k}{3^k}:x\in\{0,1+t,2-t\}\right\} .$$
Since $S_t$ is just a scaled copy of $C_{(2-t)/(1+t)}$, by Proposition 2 we have the measure
\[   \mu _1(S_t)=\left\{
\begin{array}{ll}
      (1+t)/q & (2-t)/(1+t)=p/q\text{ in lowest terms, }p+q\equiv 0\text{ (mod 3)} \\
      0 & \text{otherwise}
\end{array} 
\right. \]
Since $\mu _1(S_t)$ takes on nonzero value at only countably many $t$, the result follows from Fubini's theorem. $\hfill\Box$

\section{Parallelogram formulation}

Over 30 yeara after \cite{kaczynski67}, Kaczynski would outline a reformulation of the problem, offering three conjectures on the matter \cite{kaczynski98}. Rather than conceiving the unrestricted trapezoid as a union of lines, Kaczynski considered unions of increasingly thinner parallelograms.
\begin{definition}
Let $P_{j,k}^n=\{(x,y)\in\mathbb{R}^2:nx+(k-j)y\ge k-1,nx+(k-j)y\le k,y\in[0,1]\}$ be a parallelogram with parallel sides $[(j-1)/n,j/n]\times\{ 0\}$ and $[(k-1)/n,k/n]\times\{ 1\}$. For $\sigma\in\text{Sym}(n)$, we say that the \emph{unrestricted trapezoid} $T^n_\sigma$ is given by the union $T^n_\sigma =\bigcup _{j=1}^n P^n_{j,\sigma(j)}$.
\end{definition}
In this way, the unrestricted trapezoid associated with $f_1$ above can be expressed as $T^3_{(2~3)}$.

The following conjectures center around properties of the smallest unrestricted trapezoid of $n$ parallelograms of width $1/n$, namely $\alpha (n)=\min _{\sigma\in\text{Sym}(n)}\mu _2(T^n_\sigma)$.
\begin{conjecture}
$\lim _{n\rightarrow\infty}\alpha (n)=0$.
\end{conjecture}
Kaczynski claimed to have proved this. We will show that this follows from Propoition 1 in the following section.
\begin{conjecture}
There exists constants $c_\pm >0$ such that $\frac{c_-}{\log{n}}<\alpha (n)<\frac{c_+}{\log{n}}$.
\end{conjecture}
Kaczynski claimed to be able to prove that $c_-$ exists and also claimed a weaker upper bound in the form of $\alpha (n)<c_+\left(\log (\log n)\right)^2/\log{n}$. We will prove a different upper bound.
\begin{conjecture}
The sequence $\alpha (n)$ decreases monotonically.
\end{conjecture}
This was left as an open problem in \cite{kaczynski98}.

\section{Proving the conjectures}

To prove conjectures from the previous section, we first connect projections of the Sierpinski gasket \cite{kenyon97} to our modified Cantor sets.
\begin{definition}
We define the \emph{Sierpinski gasket} $\mathcal{G}\subset\mathbb{R}^2$ as
$$\mathcal{G}=\left\{\sum _{k=1}^\infty\frac{x_k}{3^k}:x\in\{ (0,0),(2,0),(0,2)\}\right\} .$$
We define the \emph{$n^\text{th}$ partial Sierpinski gasket} $\mathcal{G}_n\supset\mathcal{G}$ similarly. Let $G_n=\left\{\sum _{k=1}^n\frac{x_k}{3^k}:x\in\{ (0,0),(2,0),(0,2)\}\right\}$ and let $\Delta (a,b,c)$ be the area bounded by the triangle with vertices $a,b,c\in\mathbb{R}^2$.
$$\mathcal{G}_n=\bigcup_{x\in G_n}\Delta\left(x,x+\left(\frac{1}{3^n},0\right) ,x+\left( 0,\frac{1}{3^n}\right))\right)$$
\end{definition}
Let the projection operator $\text{proj}_\theta :\mathbb{R}^2\rightarrow\mathbb{R}$ map points $(x,y)\in\mathbb{R}^2$ onto a line through the origin with angle $\theta$ to the $x$-axis such that $\text{proj}_\theta :(x,y)\mapsto x\cos\theta +y\sin\theta$. We are interested in the sets $\text{proj}_\theta (\mathcal{G}_n)$; more generally there is an active field of research surrounding linear projections of fractals in $\mathbb{R}^2$ \cite{laba12}. One of the major subjects of interest is the following:
\begin{definition}
The $\emph{Favard distance}$ of a compact set $E\subset\mathbb{R}^2$ is given by
$$\text{Fav}(E)=\frac{1}{\pi}\int _0^\pi\mu\left(\text{proj}_\theta (E)\right)d\theta$$
\end{definition}
While it is perhaps obvious that $\text{Fav}(\mathcal{G})=0$ from proposition 2, a result by Bond and Volberg \cite{bond09} gives a decay bound on the Favard distances of partial Sierpinski gaskets.
\begin{proposition}
There exist constants $C,p>0$ such that $\text{Fav}(\mathcal{G}_n)\le Cn^{-p}$.
\end{proposition}
A later result by Bond and Volberg proves that $p>1/14$ \cite{bond10}.

The rest of this section centers around utilizing proposition 3 to prove the conjectures from section 3. To do this, we must compare $\text{proj}_\theta (\mathcal{G}_n)$ to the modified Cantor sets from Proposition 2 and the horizontal slices of the unrestricted trapezoid. Let us define partial variants of these sets:
\begin{definition}
Let $D_n(a,b,c)=\left\{\sum _{k=1}^n\frac{x_k}{3^k}:x_k\in\{ a,b,c\}\right\}$. If $T^{3^n}_\sigma$ is the unrestricted trapezoid associated with the $n^\text{th}$ iteration induced by swapping $1\leftrightarrow 2$ in a base 3 expansion (i.e. if $x-1=x_nx_{n-1}...x_1$ is an $n$-digit base 3 number, then $\sigma (x-1)+1=y_ny_{n-1}...y_1$ where $y_k=2x_k\text{ (mod 3)}$), then we write $S^{(n)}_t=(\mathbb{R}\times\{ t\} )\cap T^{3^n}_\sigma$, the slice of $T^{3^n}_\sigma$ at $y=t$, as
$$S^{(n)}_t=\bigcup _{x\in D_n(0,1+t,2-t)}\left[ x,x+\frac{1}{3^n}\right]$$
We similarly define the \emph{$n^\text{th}$ partial modified Cantor set} as
$$C^{(n)}_t=\bigcup _{x\in D_n(0,1,t)}\left[ x,x+\frac{1}{3^n}\right]$$
\end{definition}
This allows us to make the following comparison:
\begin{lemma}
For $t\in [0,1]$ and $\varphi (t)=\tan^{-1}\left(\frac{2-t}{1+t}\right)$, we have
$$\mu (S^{(n)}_t)\le (1+t)\mu\left(\text{proj}_{\varphi (t)}(\mathcal{G}_n)\right)$$
\end{lemma}
{\bf Proof:} We can show that $\mu (S^{(n)}_t)\le (1+t)\mu\left( C^{(n)}_{(2-t)/(1+t)}\right)$ for $t\ge 0$ since
\begin{align*}
S^{(n)}_t &= \bigcup _{x\in D_n(0,1+t,2-t)}\left[ x,x+\frac{1}{3^n}\right]\subset\bigcup _{x\in D_n(0,1+t,2-t)}\left[ x,x+\frac{1+t}{3^n}\right] \\
    &= \bigcup _{x\in D_n(0,1,\frac{2-t}{1+t})}(1+t)\left[ x,x+\frac{1}{3^n}\right] =(1+t)C^{(n)}_{(2-t)/(1+t)}
\end{align*}
The result then follows from noticing that for $\theta\le\frac{\pi}{4}$, we have $\text{proj}_\theta (\mathcal{G}_n)=C^{(n)}_{\tan\theta}$. $\hfill\Box$

Combining these results leads us to a similar polynomial bound on the decay of the unrestricted trapezoids induced by the iterative scheme described in Definition 6.
\begin{proposition}
Let $\sigma$ be the permutation on base 3 numbers as defined in Definition 6. Then there exist constants $C,p>0$ such that $\mu\left( T^{3^n}_\sigma\right)\le Cn^{-p}$.
\end{proposition}
{\bf Proof:} By Fubini's theorem, The measure of $T^{3^n}_\sigma$ is given by integrating over the measures of the horizontal slices
$$\mu _2\left( T^{3^n}_\sigma\right) =\int _0^1\mu\left( S^{(n)}_t\right) dt$$
From Lemma 1, this integral can be bound by an integral over projections of the Sierpinski gasket:
$$\int _0^1\mu\left( S^{(n)}_t\right) dt\le \int _0^1(1+t)\mu\left(\text{proj}_{\varphi (t)}(\mathcal{G}_n)\right) dt$$
By the change of variables $\theta =\varphi (t)$, we can write:
$$\int _0^1(1+t)\mu\left(\text{proj}_{\varphi (t)}(\mathcal{G}_n)\right) dt =\int _{\tan ^{-1}(1/2)}^{\tan ^{-1}(2)}\frac{9(\tan ^2\theta +1)}{\left(\tan\theta +1\right) ^3}\mu\left(\text{proj}_{\theta}(\mathcal{G}_n)\right) d\theta$$
Placing an estimate on the trigonometric portion of the integral we can write
$$\int _{\tan ^{-1}(1/2)}^{\tan ^{-1}(2)}\frac{9(\tan ^2\theta +1)}{\left(\tan\theta +1\right) ^3}\mu\left(\text{proj}_{\theta}(\mathcal{G}_n)\right) d\theta\le\frac{10}{3}\int _{\tan ^{-1}(1/2)}^{\tan ^{-1}(2)}\mu\left(\text{proj}_{\theta}(\mathcal{G}_n)\right) d\theta$$
Since Lebesgue measure is nonnegative, we can extend the limits of integration to bound the measure of $T^{3^n}_\sigma$ by a factor of the associated Favard distance
$$\mu\left( T^{3^n}_\sigma\right)\le\frac{10\pi}{3}\text{Fav}\left(\mathcal{G}_n\right)$$
Proposition 3 completes the proof. $\hfill\Box$

This result translates to a logarithmic bound on $\alpha (n)$ \emph{only} for $n=3^m$ where $m\in\mathbb{N}$. If $n$ has the base 3 representation $n=x_kx_{k-1}...x_0$, we can define $\sigma _n$ as in Figure 3 such that the first $x_k3^k$ intervals map to $x_k$ copies of $T^{3^k}_\sigma$, the next $x_{k-1}3^{k-1}$ intervals map to $x_{k-1}$ copies of $T^{3^{k-1}}_\sigma$, etc. Using this scheme, we can prove a weaker upper bound on Conjecture 2 (consequently proving Conjecture 1 as well):
\begin{theorem}
There exist constants $C,p>0$ such that $\alpha (n)<C/(\log{n} )^{p}$.
\end{theorem}
Before proving this theorem, we supply a technical lemma:
\begin{lemma}
For $p\in (0,1)$, there exists a constant $C$ such that
$$\int _0^ne^xx^{-p}dx\le Ce^nn^{-p}$$
\end{lemma}
{\bf Proof:} Let us make the change of variables $t=1-x/n$ such that this integral can be transformed into
$$\int _0^ne^xx^{-p}dx=\frac{e^n}{n^{p-1}}\int _0^1\frac{e^{-nt}}{(1-t)^p}dt$$
For any constant $C>1$, we can find $t_0$ dependent only on $C,p$ such that $pt+1>(1-t)^{-p}$ for $t\in (0,t_0)$. This allows us to split the integral
$$\int _0^1\frac{e^{-nt}}{(1-t)^p}dt\le \int _0^{t_0}C(pt+1)e^{-nt}dt+\int _{t_0}^{1}\frac{e^{-nt_0}}{(1-t)^p}dt$$
Both of these integrals have closed form expressions, the leading asymptotic term $1/n$ arising from the left integral, thus completing the proof. $\hfill\Box$

\begin{figure}
\centering
\includegraphics[scale=0.2]{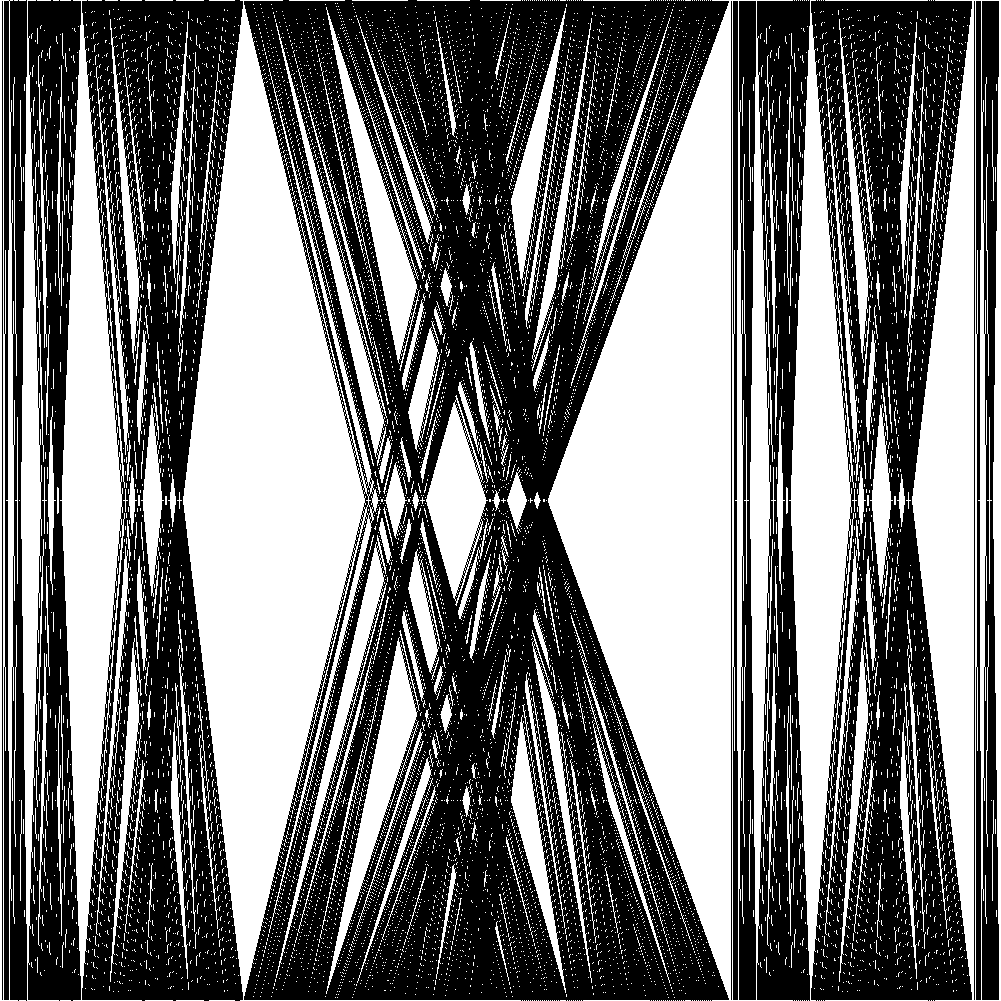}\hspace{1cm}\includegraphics[scale=0.2]{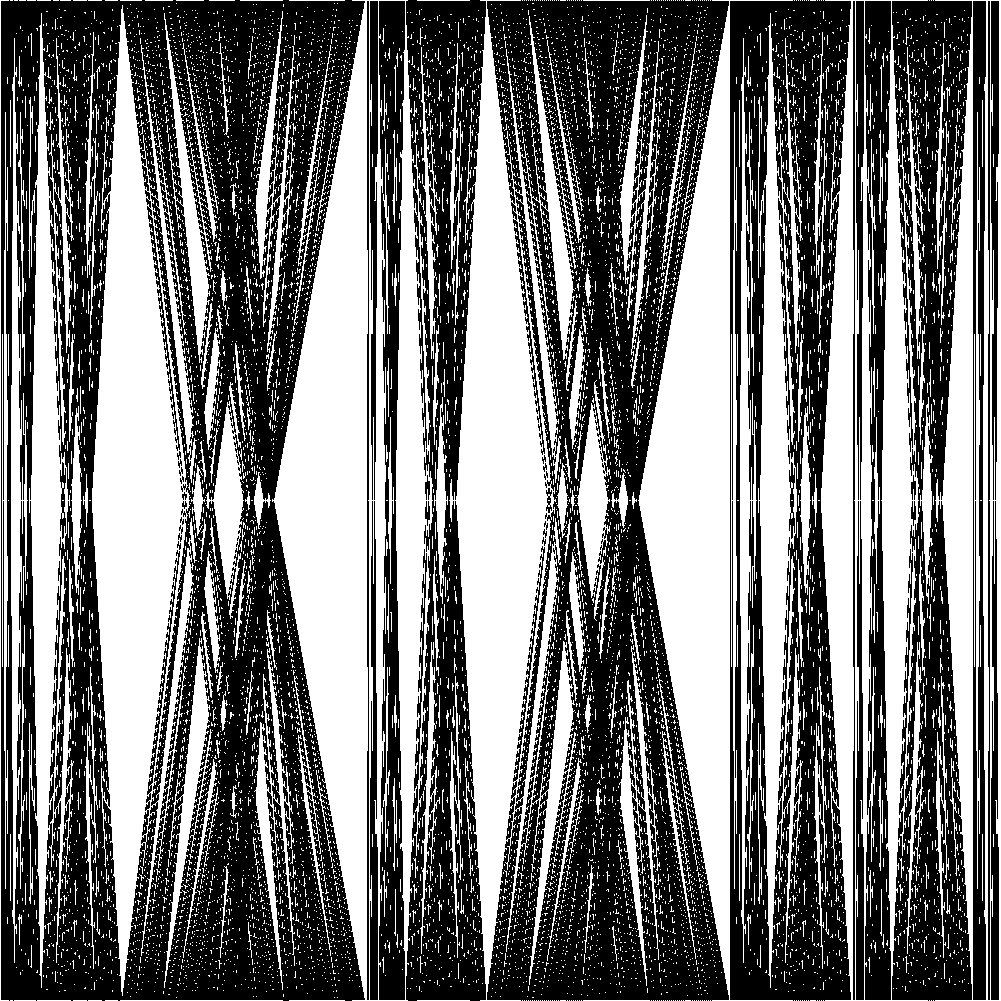}
\caption{Examples of parallelogram-defined unrestricted trapezoids $T^n_\sigma$ where $\sigma$ is as outlined in the proof of Theorem 4 (\emph{Left}) $n=1000$ which has base-3 representation 1101001 (\emph{Right}) $n=2000$ which has base-3 representation 2202002.}
\end{figure}

\noindent {\bf Proof of Theorem 4:} It is clear that $\alpha (n)\le\mu\left( T^n_{\sigma _n}\right)$. Let $n=x_kx_{k-1}...x_0$ be a base 3 integer representation. Then the following holds:
$$\mu\left( T^n_{\sigma _n}\right) =\frac{1}{n}\sum _{j=0}^k x_j3^j \mu\left( T^{3^j}_{\sigma}\right)$$
Each digit must satisfy $x_j\le 2$, and by Proposition 4 there exist constants $C,p>0$ such that $\mu\left( T^{3^j}_{\sigma}\right) \le Cj^{-p}$, so we may write
$$\frac{1}{n}\sum _{j=0}^k x_j3^j \mu\left( T^{3^j}_{\sigma}\right) <\frac{C}{n}\sum _{j=1}^{\lfloor\log _3n\rfloor} 3^jj^{-p}$$
Since the function $f(x)=3^xx^{-p}$ is increasing for $x>p/\log _3(x)$, we can bound this sum by the integral
$$\frac{C}{n}\sum _{j=1}^{\lfloor\log _3n\rfloor} 3^jj^{-p}\le\frac{C}{n}\int _0^{\log _3 n+1}3^xx^{-p}dx$$
By making the change of variables $x=u\log _3e$, we can write
$$\frac{C}{n}\int _0^{\log _3 n+1}3^xx^{-p}dx=\frac{C(\log _3e)^{1-p}}{n}\int _0^{\log 3n}e^uu^{-p}du$$
The result follows from applying Lemma 2 and collecting constants. $\hfill\Box$

\section{Application to Cluster Sets}

While Theorem 1 assisted in constructing a measurable half-plane function with nonmeasurable boundary function in \cite{kaczynski67}, the open questions from \cite{kaczynski67} and \cite{kaczynski98} answered by Propositions 1 and 4 were also initially conceived to address problems related to boundary functions and cluster sets. We state the relevant material here for completeness.
\begin{definition}
Given a function defined on the open unit disk $f:\mathbb{D}\rightarrow\mathbb{C}$ and a Jordan arc $\gamma$ with an endpoint on the boundary of the disk, we say that $C(f,\gamma )$ is the \emph{cluster set} of $f$ on $\gamma$ where
$$C(f,\gamma )=\{ w\in\mathbb{C}:\exists \{z_n\}~s.t.~z_n\in\gamma, |\lim z_n|=1, \lim f(z_n)=w\}$$
\end{definition}
\begin{definition}
Let $\gamma _1,\gamma_2,\gamma _3\subset\mathbb{D}$ be Jordan arcs with a similar endpoint $z$ on the unit circle. If, for a disk function $f$, the intersection satisfies $C(f,\gamma _1)\cap C(f,\gamma _2)\cap C(f,\gamma _3)=\emptyset$, then $f$ has the \emph{three-arc property} at $z$. If each $\gamma _1,\gamma _2,\gamma _3$ is a union of rectilinear segments, then we say $f$ has the \emph{three-segment property} at $z$.
\end{definition}

In \cite{kaczynski67}, a version of conjecture 1 was suggested to assist in answering the following open question from \cite{bagemihl59} (presented as conjecture):
\begin{conjecture}
There exists a continuous function in $\mathbb{D}$ having the three-segment property at each point of a set of positive measure or second category on $|z|=1$.
\end{conjecture}
Weaker versions of this conjecture have been proved. Jarn\'{i}k \cite{jarnik36} gave an example of a function having the three-segment property at uncountably many points. Piranian was cited in \cite{bagemihl59} as having proved the existence of a continuous disk function with the three-arc (not segment) property. Bagemihl later proved the existence of a normal meromorphic (not continuous) disk function with the three-segment property \cite{bagemihl74}.

\bibliographystyle{unsrt}
\bibliography{example}

\end{document}